\documentclass[11pt]{article}

\usepackage{amsthm}
\usepackage{latexsym}
\usepackage{amssymb}
\usepackage{amsmath}
\usepackage{mathrsfs}
\usepackage{epsfig}
\usepackage{euscript}

\newcommand{\half}{\frac{\scriptstyle 1}{\scriptstyle 2}}

\newcommand{\R}{\mathbb{R}}
\newcommand{\N}{\mathbb{N}}

\renewcommand{\d}{\mathrm{d}}

\newcommand{\id}{\mathrm{Id}}
\newcommand{\supp}{\mathrm{supp}}

\newcommand{\tx}{\tilde{X}}
\newcommand{\ResSig}{_{|_\Sigma}}
\newcommand{\Tsig}{\mathbb{T}_{_\Sigma}}
\newcommand{\dnu}{\mathrm{d} \nu}
\newcommand{\dsig}{\mathrm{d} \nu_{_\Sigma}}
\newcommand{\chardsig}{\mathrm{d} \nu^0_{_\Sigma}}

\newcommand{\ih}{\,^i \hskip-0.15em h}
\newcommand{\gk}{\,^k \hskip-0.18em g}
\newcommand{\oplk}{{\,^k \hskip-0.20em L_1}}
\newcommand{\bk}{\,^k \hskip-0.12em b}
\newcommand{\ck}{\,^k \hskip-0.15em c}
\newcommand{\lambdak}{\,^k \hskip-0.15em \lambda}

\newtheorem{theorem}{Theorem}
\newtheorem{proposition}{Proposition}[section]
\newtheorem{corollary}{Corollary}[section]

\newtheorem{remark}{Remark}[section]
\newtheorem{property}{Property}

\begin{document}
\begin{center}

{\huge{\bf On Lars H\"ormander's remark on the}}

\vspace{0.1in}

{\huge{\bf  characteristic Cauchy problem}}

\vspace{0.1in}

{\large{\bf Jean-Philippe NICOLAS}}
\\
{\small{\it{Institut de Math\'ematiques, M.A.B.,
Universit\'e Bordeaux 1,}}} \\
{\small{\it{351 Cours de la Lib\'eration, 
33405 Talence Cedex, France. }}} \\
{\small{\it{Jean-Philippe.Nicolas@math.u-bordeaux1.fr}}}
\end{center}

\vspace{0.1in}
\begin{center}
{\bf Abstract}
\end{center}

We extend the results of a work by L. H\"ormander
  \cite{Ho} concerning the resolution of the characteristic Cauchy
  problem for second order wave equations with regular first order
  potentials. The geometrical background of this work was a spatially
  compact spacetime with smooth metric. The initial data surface was
  spacelike or null at each point and merely Lipschitz. We lower the
  regularity hypotheses on the metric and potential and obtain similar
  results. The Cauchy problem for a spacelike initial data surface is
  solved for a Lipschitz metric and coefficients of the first order
  potential that are $L^\infty_\mathrm{loc}$, with the same finite
  energy solution space as in the smooth case. We also solve the fully
  characteristic Cauchy problem with very slightly more regular metric
  and potential~: essentially, a ${\cal C}^1$ metric and a potential
  with continuous coefficients of the first order terms and locally $L^\infty$
  coefficients for the terms of order $0$.

\begin{center}
{\bf R\'esum\'e}
\end{center}

Nous \'etendons des r\'esultats dus \`a L. H\"ormander \cite{Ho}
concernant la r\'esolution du probl\`eme de Cauchy caract\'eristique
pour des \'equations d'onde du second ordre avec un potentiel
r\'egulier du premier ordre. Le cadre g\'eom\'etrique de \cite{Ho}
\'etait un espace-temps spatialement compact avec une m\'etrique
r\'eguli\`ere. L'hypersurface sur laquelle les donn\'ees initiales
sont fix\'ees \'etait spatiale ou caract\'eristique en chaque point et
simplement de r\'egularit\'e Lipschitz. Nous affaiblissons les
hypoth\`eses de r\'egularit\'e sur la m\'etrique et le potentiel et
nous obtenons des r\'esultats analogues. Le probl\`eme de Cauchy pour
une hypersurface spatiale est r\'esolu dans le cas d'une m\'etrique
Lipschitz et pour un potentiel dont les coefficients sont localement
$L^\infty$, avec le m\^eme espace de solutions que dans le cas
r\'egulier. Nous r\'esolvons \'egalement le probl\`eme de Cauchy
totalement caract\'eristique dans un cadre tr\`es l\'eg\`erement plus
r\'egulier~: essentiellement, une m\'etrique ${\cal C}^1$ et un
potentiel dont les coefficients des termes du premier ordre sont
continus et ceux des termes d'ordre $0$ sont localement $L^\infty$.

\section{Introduction}

In 1990, in a paper entitled ``A remark on the characteristic Cauchy
problem'' \cite{Ho}, L. H\"ormander solved the global Cauchy problem
for a class of wave equations on spatially compact space-times with
initial data hypersurfaces that were weakly spacelike (i.e. at almost
every point either null or spacelike) and merely Lipschitz. He assumed
the metric on his space-time to be smooth and the equations he was
considering were first order perturbations of the covariant scalar
wave equation associated with the metric, the first order perturbation
consisting of a differential operator with smooth coefficients. At the
end of his work, he remarked that all the estimates depended only on
the Lipschitz norm of the metric and the $L^\infty$ norms of the
coefficients of the first order terms (on a given large enough time
interval on which the solution is studied)~; he concluded that this
was the proper generality of his theorem. However, he did not prove
that the result could be extended to the case of a Lipschitz metric
and a first order differential perturbation with $L^\infty$
coefficients. To this day and to the author's knowledge, this final
remark has remained unchecked.

In addition to the mathematical interest of this question, it is
relevant for matters related to scattering theory in general
relativity in which the author is at present involved~: namely, the
construction of geometrical versions of scattering theory in generic
non-stationary space-times, based on conformal techniques. Such ideas
can be traced back to 1963, when R. Penrose put forward in \cite{Pe}
the essential ideas of the resolution of the characteristic Cauchy
problem for field equations in relativity. These ideas were then used
by F.G. Friedlander to give a conformal construction of scattering
theory for the wave equation on static space-times in
\cite{Fri1,Fri2}, respectively in 1980 and 2001. In 1990, J.C. Baez,
I.E. Segal and Z.F. Zhou applied a similar idea to nonlinear
conformally invariant wave equations on flat space-time. On the same
year and in the same issue of the same journal, L. H\"ormander
produced his rigorous proof of the solution to the characteristic Cauchy
problem for wave equations on spatially compact space-times
\cite{Ho}. This is fundamental for any conformal description of
scattering and H\"ormander's proof, based on energy estimates, allows
to work with minimum regularity solutions, which is the natural
framework of scattering theory. In 2004, L.J. Mason and the author
\cite{MaNi} proposed a conformal construction of a scattering operator
for spin 1/2 and spin 1 massless fields on generic non stationary
asymptotically simple space-times with smooth conformal infinity. This
construction uses Penrose's ideas and a modification of H\"ormander's
proof that allows to establish the equivalence with a standard analytic
time-dependent scattering theory, defined in terms of classical wave
operators.

The notion of asymptotic simplicity, on which \cite{MaNi} was strongly
relating, was defined by R. Penrose, as a simple geometric description
of asymptotic flatness. For some time, this was considered by many as
an idealised model, because no known solution of the Einstein vacuum
equations, except Minkowski space-time, was asymptotically simple. The
first examples of vacuum space-times that approached the
asymptotically simple framework were the space-times of
D. Christodoulou and S. Klainerman \cite{ChriKla}, propagated from
initial data close to Minkowski space. These space-times are
non-stationnary and only fail to be asymptotically simple because of a
small lack or regularity at null infinity. Since this work, other
generically non-stationary vacuum space-times have been constructed,
by J. Corvino \cite{Co}, P. Chrusciel and E. Delay
\cite{ChruDe1,ChruDe2} and J. Corvino and R.M. Schoen
\cite{CoScho}. They are asymptotically simple with specifyable
regularity at null and timelike infinity, and diffeomorphic to
Schwarzschild or Kerr spacetime in a neighbourhood of spacelike
infinity. Then, S. Klainerman and F. Nicol\`o \cite{KlaNi} proved that
for initial data that are close to flat space-time and with
stronger fall-off assumptions than in \cite{ChriKla}, the
corresponding global solution of the Einstein vacuum equations is
asymptotically simple.

The regularity of conformal infinity encodes the fall-off of
the physical metric~; the more regular conformal infinity, the
stronger the fall-off. So, working with less regularity of the
conformal metric at conformal infinity, means working with larger
classes of physical metrics. If in particular one wishes to extend the
constructions of \cite{MaNi} to the space-times of \cite{ChriKla},
then the characteristic Cauchy problem must be solved in frameworks
that are only slightly more regular than what L. H\"ormander had in
mind in his final remark in \cite{Ho}. It is therefore quite crucial
to prove this remark. The present work is a step in this
direction. The results of \cite{Ho} are extended to the regularity
setting proposed by H\"ormander for the standard Cauchy problem~; we
also treat the fully characteristic Cauchy problem (Goursat problem),
for this we impose a little more regularity, but we remain below the
actual regularity of conformal infinity for the space-times of
\cite{ChriKla}. The paper is organized as follows.
\begin{itemize}
\item In section \ref{SectionHormRes}, we describe the geometrical setting
      and give a summary of L. H\"ormander's result (theorem
      \ref{thmHorm}). His work deals
      with the Cauchy problem for initial data surfaces that are
      allowed locally to be spacelike or null and thus solves the
      standard Cauchy problem as well as the characteristic Cauchy
      problem.
\item Section \ref{MainResults} contains the main results of this
      work. In subsection \ref{SectionCauchy}, theorem \ref{Cauchy}
      solves the Cauchy problem in the regularity setting proposed by
      L. H\"ormander. The surprising thing is that the minimum
      regularity solutions in fact remain continuous in time with
      values in $H^1$ instead of being only locally $L^\infty$ with
      values in $H^1$. A corollary for equations that are homogeneous
      of the second order is then obtained. In this case, we have
      access to $H^2$ solutions for more regular data.
      Subsection \ref{SectionGoursat} treats the fully characteristic Cauchy
      problem. It uses crucially the corollary of section
      \ref{SectionCauchy} to define a regularization of solutions that
      has the strong convergence properties appropriate for obtaining the
      fundamental energy estimates. In theorem \ref{Goursat}, the well
      posedness is proved for a ${\cal C}^1$ metric and coefficients
      of lower order terms that
     are assumed to be continuous for the first order terms and
     $L^\infty_\mathrm{loc}$ for the zero-order terms. The proof
     follows the essential structure of H\"ormander's proof, based on
     two reciprocal energy estimates between some spacelike slice
     and the null hypersurface, followed by the construction of one
     solution to the characteristic Cauchy problem. It turns out that
     the theorem is in fact valid for a regularity setting
     intermediate between what H\"ormander proposed and that of
     theorem \ref{Goursat}~; this is expressed in theorem
     \ref{TheorMinReg}.
\item The proofs of the theorems are given in the last section.
\end{itemize}

For simplicity, we work with a scalar wave equation with real-valued
unknown function. However, the theorems are also valid for a wave
equation with complex, tensor or spinor valued unknown function (for
spinor fields, provided the space-time admits a spin-structure).

\section{Geometrical and functional framework and summary of
  Lars H\"ormander's result}
\label{SectionHormRes}

The geometrical framework chosen by H\"ormander is as follows~: $X$ is
a ${\cal C}^\infty$ compact manifold of dimension $n\geq 1$ and $\tx =
\R_t \times X$. For $t\in \R$, we denote $X_t = \{ t \} \times X$. We
consider on $X$ a time dependent Riemannian metric $g(t)$ assumed to
be ${\cal C}^\infty$ on $\tx$. An immediate consequence of this is
\begin{property} \label{ControlMetric}
There exist two continuous positive functions\footnote{In \cite{Ho},
  for simplicity, it is assumed that $C_1$ and $C_2$ are
  constants. This is of course unimportant since one can always
  restrict the study to a generic compact time interval $[-T , T]$.}
$C_1$ and $C_2$ on $\R$ such that, for a given local smooth coordinate
system on $X$, the matrix $G=\left( g_{\alpha \beta} \right)$
satisfies, as a quadratic form on $\R^n$~:
\[ C_1 (t) \id_n \leq G(t,x) \leq C_2
(t) \id_n \, ,~~ \forall (t,x) \in \tx \, , \]
where $\id_n$ denotes the $n\times n$ identity matrix.
\end{property}
We also define $\d \nu$ a fixed smooth density on $X$~; in local
coordinates $\dnu = \gamma \d x$. We can assume that $\dnu$ is the
volume measure induced by a smooth Riemannian metric $h$ on $X$~; $\d
t \dnu$ is then the measure induced by the smooth Riemannian metric
$\tilde{h} = \d t^2 + h$ on $\tx$. We denote $\nabla$ and
$\tilde{\nabla}$ the Levi-Civita connections induced respectively by
$h$ and $\tilde{h}$.

We work with Sobolev spaces $H^\mu$ and $H^\mu_\mathrm{loc}$ defined
on $X$ and $\tx$ for any $\mu \in \R$ by local identification with the
corresponding function spaces on smooth open sets of $\R^n$ and
$\R^{n+1}$. We only use explicit norms for $\mu = 0$ or $1$~; in
fact we have natural norms on $H^k(X)$ and $H^k(\tx )$ for any $k\in
\N$~:
\begin{equation} \label{HkXNorm}
\| u \|^2_{H^k(X)} = \sum_{p=0}^k \int_X \left< \nabla^p u ,
    \nabla^p u \right> \dnu \, ,
\end{equation}
\begin{equation} \label{HktxNorm}
\| u \|^2_{H^k(\tx)} = \sum_{p=0}^k \int_{\tx} \left<
    \tilde{\nabla}^p u , \tilde{\nabla}^p u \right> \d t \dnu
\end{equation}
where the same notation $<.,.>$ refers to the inner product on tensors
at a point induced by $h$ or $\tilde{h}$. The inner products
associated to the norms (\ref{HkXNorm}) and (\ref{HktxNorm}) are
denoted $<.,.>_{H^k (X)}$ and $<.,.>_{H^k (\tx)}$. On $H^1 (X)$, we
also define a norm $\| . \|_{H^1 (X_t )}$ that is more
closely related to the metric $g(t)$~:
\begin{equation} \label{H1PhysNorm}
\| u \|^2_{H^1 (X_t )} = \int_{X} \left( g^{\alpha \beta} (t,x)
    \partial_\alpha u (x) \partial_\beta u (x) + \left| u(x)
    \right|^2 \right) \dnu (x) \, .
\end{equation}
The $H^1$ norms (\ref{HkXNorm}) and (\ref{H1PhysNorm}) are equivalent
for any $t\in \R$ and the equivalence is locally uniform in $t$.
Another type of function space we shall need to consider is
$W^{1,\infty} ( {\cal O})$ (resp. $W^{1,\infty}_\mathrm{loc} ( {\cal
  O})$), where $\cal O$ is an open set of $X$ or $\tx$~; it is defined
as the space of functions in $L^\infty ({\cal O})$
(resp. $L^\infty_\mathrm{loc} ({\cal O})$) such that their gradient is
also in $L^\infty ({\cal O})$ (resp. $L^\infty_\mathrm{loc} ({\cal
  O})$).

On $\tx$, we consider a wave equation of the form
\begin{equation} \label{WaveEq}
\square u + L_1 u = 0
\end{equation}
where $\square$ denotes the simplified d'Alembertian
\begin{equation} \label{Box}
\square = \frac{\partial^2}{\partial t^2} - \gamma^{-1}
\frac{\partial}{\partial x^\alpha} \left( \gamma g^{\alpha \beta}
  \frac{\partial}{\partial x^\beta} \right) \, ,
\end{equation}
and $L_1$ is a general first order differential operator
\begin{equation} \label{L1}
L_1 = b^0
\frac{\partial}{\partial t} + b^\alpha \frac{\partial}{\partial
  x^\alpha} + c
\end{equation}
whose coefficients $b^0$, $b^\alpha$ and $c$ are assumed to be ${\cal
  C}^\infty$ functions on $\tx$. The hypersurface on which the initial
data are specified can be a spacelike Cauchy hypersurface for a
standard Cauchy problem, a light cone for a characteristic Cauchy
problem (Goursat problem), or anything in between. It is defined as
follows
\begin{equation} \label{Sigma}
\Sigma = \left\{ (\varphi (x ) , x) \, ;~x\in X \right\} \, ,~\varphi
~:~X\longrightarrow \R \, ,
\end{equation}
where $\varphi$ is simply assumed to be Lipschitz on $X$, to allow for
singularities such as the vertex of a light cone, and weakly spacelike,
i.e.
\begin{equation} \label{WeakSpacelike}
g^{\alpha \beta} (\varphi (x) , x) \partial_\alpha \varphi (x)
\partial_\beta \varphi (x) \leq 1 ~\mathrm{almost~everywhere~on~}X \,
.
\end{equation}
Condition (\ref{WeakSpacelike}) has a meaning, since Lipschitz
functions are differentiable almost everywhere, and it simply says
that $\Sigma$ is allowed to be locally spacelike or null but not
timelike.

We consider on $\Sigma$ the density measure $\dsig$ which is simply
$\d \nu$ lifted to $\Sigma$ using parametrization (\ref{Sigma}). 
The hypersurface $\Sigma$ being merely Lipschitz, we can define
the spaces $H^\mu (\Sigma )$ only for $|\mu | \leq 1$~; these spaces
are canonically isomorphic to the corresponding Sobolev spaces on $X$
by (\ref{Sigma}). On $L^2 (\Sigma )$ and $H^1
(\Sigma )$, we consider the norms $\| . \|_{L^2 (\Sigma )}$ and $\|
. \|_{H^1 (\Sigma )}$, naturally induced by this isomorphism. We also
define a norm $\left\| . \right\|_{H^1 (\Sigma ; g)}$on $H^1 (\Sigma
)$ in the two following equivalent manners~: if the element $\psi$ of
$H^1 (\Sigma )$ is considered as the lift on $\Sigma$ of an element of
$H^1 (X)$, the norm has the form
\begin{equation} \label{PhysNorm1}
\left\| \psi \right\|^2_{H^1 (\Sigma ; g)} = \int_X \left\{ \left|
    \psi \right|^2 + g^{\alpha \beta} (\varphi (x),x ) \partial_\alpha
    \psi (x) \partial_\beta \psi (x) \right\} \dnu
\end{equation}
and if $\psi$ is defined as the trace on $\Sigma$ of some $\Psi \in
H^{3/2}_\mathrm{loc} (\tx)$,
\begin{equation} \label{PhysNorm2}
\left\| \psi \right\|^2_{H^1 (\Sigma ; g)} = \int_\Sigma \left\{
    \left| \Psi \right|^2 + g^{\alpha \beta} \left( \partial_\alpha
    \Psi + \partial_\alpha \varphi \partial_t \Psi \right)  \left(
    \partial_\beta \Psi + \partial_\beta \varphi \partial_t
    \Psi \right) \right\} \dsig \, .
\end{equation}
The norms $\left\| . \right\|_{H^1 (\Sigma ; g)}$ and $\left\|
  . \right\|_{H^1 (\Sigma )}$ are of course equivalent. We shall also
consider the foliation $\{ \Sigma_t \}_{t\in \R}$
\begin{equation} \label{Sigmat}
\Sigma_t = \left\{ \left( t+\varphi (x) , x \right) \, ; x\in X
\right\} \, , ~\Sigma_0 = \Sigma \, .
\end{equation}
On each $\Sigma_t$, we define the spaces $H^\mu (\Sigma_t )$, $-1 \leq
\mu \leq 1$. These spaces are canonically isomorphic to the
corresponding spaces on $\Sigma$ by parametrizations (\ref{Sigma}) and
(\ref{Sigmat}). We use this canonical isomorphism to identify $H^\mu
(\Sigma_t )$ with $H^\mu (\Sigma )$. On $H^1 (\Sigma_t )$, in addition
to the norm $\| . \|_{H^1 (\Sigma )}$ inherited from the previous
identification, we can also consider a norm involving the restriction
to $\Sigma_t$ of the metric $g$. Its definition is analogous to
(\ref{PhysNorm1}) and (\ref{PhysNorm2})~: for $\psi \in H^1 (\Sigma_t
)$ seen as the lift on $\Sigma_t$ of an element of $H^1 (X)$,
\begin{equation}
\left\| \psi \right\|^2_{H^1 (\Sigma_t ; g)} = \int_X \left\{ \left|
    \psi \right|^2 + g^{\alpha \beta} \left( t+\varphi (x),x \right)
    \partial_\alpha \psi (x) \partial_\beta \psi (x)
    \right\} \dnu
\end{equation}
and for $\psi$ defined as the trace on $\Sigma_t$ of some $\Psi \in
H^{3/2}_\mathrm{loc} (\tx)$,
\begin{equation}
\left\| \psi \right\|^2_{H^1 (\Sigma_t ; g)} = \int_{\Sigma_t} \left\{
    \left| \Psi \right|^2 + g^{\alpha \beta} \left( \partial_\alpha
    \Psi + \partial_\alpha \varphi \partial_t \Psi \right)  \left(
    \partial_\beta \Psi + \partial_\beta \varphi \partial_t
    \Psi \right) \right\} \dsig \, .
\end{equation}
These norms are equivalent, locally uniformly in time, with the $H^1
(\Sigma )$ norm.

The well-posedness of the Cauchy problem for (\ref{WaveEq}) in $H^1(X)
\oplus L^2 (X)$ is well-known~: for any initial data $(u_0 ,u_1 ) \in
H^1(X) \oplus L^2 (X)$, for any initial time $s\in \R$, (\ref{WaveEq})
admits a unique solution $u$ in
\begin{equation}
{\cal F} = {\cal C}^0 \left( \R_t \, ; H^1 (X) \right) \cap {\cal C}^1
\left( \R_t \, ; L^2 (X) \right)
\end{equation}
such that $u(s) = u_0$ and $\partial_t u (s) = u_1$. For $u \in {\cal
  F}$, we introduce the energy of $u$ at time $t$ as the norm of
  $(u(t) , \partial_t u(t) )$ in $H^1 \oplus L^2$~:
\begin{equation} \label{Energy}
E(t,u) = \left\| u(t) \right\|_{H^1 (X_t)}^2 + \left\| \partial_t u(t)
\right\|_{L^2 (X)}^2 = \int_{X_t} \left\{ \left| \partial_t u
\right|^2 + g^{\alpha \beta} \partial_\alpha u \partial_\beta \bar{u}
+ \left| u \right|^2 \right\} \dnu \, .
\end{equation}
If $u \in {\cal F}$ is a solution of (\ref{WaveEq}), it satisfies for
all $T>0$ the energy estimate
\begin{equation} \label{EnEst1}
E(t,u) \leq E(s,u) e^{K_1 (T,g,L_1) |t-s|} ~\forall t,s \in [-T , T]
\end{equation}
where $K_1$ is a continuous positive function of $T>0$, the norms in
$W^{1,\infty} (]-T,T[ \times X )$ of $g$ and $g^{-1}$ and the norms of
the coefficients of $L_1$ in $L^{\infty} (]-T,T[ \times X )$. We
denote by $\cal E$ the space of finite energy solutions of
(\ref{WaveEq}), i.e. the set of solutions of (\ref{WaveEq}) in $\cal
F$. The energy estimate (\ref{EnEst1}) shows that for any $t \in \R$
and for any $T>0$, the following are equivalent norms on $\cal E$~:
\begin{equation} \label{EnNorm}
N(t) \, : ~u \in {\cal E} \longmapsto \sqrt{E(t,u)}
\end{equation}
and
\begin{equation} \label{UnifEnNorm}
\| u \|_{{\cal F} , T} := \sup_{-T <\tau <T} N(\tau )(u) \, .
\end{equation}

The main result of \cite{Ho} is the following~:
\begin{theorem} {\bf (H\"ormander, 1990)} \label{thmHorm}
We define on $\Sigma$ the density measure
\[ \chardsig = \left( 1 - g^{\alpha \beta} \partial_\alpha \varphi
  \partial_\beta \varphi \right) \dsig \, ,\]
which is positive where $\Sigma$ is spacelike and vanishes where
$\Sigma$ is null, and the associated $L^2$ space $L^2 (\Sigma ;
\chardsig )$. The application
\begin{equation} \label{IsomSmooth}
\begin{array}{cccc}
{\Tsig ~: }& {\cal E} & \longrightarrow & {H^1 (\Sigma ) \oplus L^2
 (\Sigma ; \chardsig )} \\
 & u & \longmapsto & {\left( u\ResSig ~ ,~ \partial_t u \ResSig
   \right) \, ,} \end{array}
\end{equation}
which is well defined for smooth solutions, extends as an
isomorphism. In particular, there exist $K_2 (T,g,L_1 )$ and $K_3
(T,g,L_1 )$, two positive continuous functions of $T>0$, the norms in
$W^{1,\infty} (]-T,T[ \times X )$ of $g$ and $g^{-1}$ and the norms of
the coefficients of $L_1$ in $L^{\infty} (]-T,T[ \times X )$, such that
for $u\in {\cal E}$, for $T>0$ satisfying $-T < \min \{ \varphi (x) \,
,~x \in X \}$, $T > \max \{ \varphi (x) \, ,~x \in X \}$, we have
\begin{equation} \label{EnEst2}
\left\| \Tsig u \right\|_{1,\Sigma} \leq K_2 (T,g,L_1 ) \left\| u
\right\|_{{\cal F},T}
\end{equation}
and
\begin{equation} \label{EnEst3}
 \left\| u \right\|_{{\cal F},T} \leq K_3 (T,g,L_1 ) \left\| \Tsig
 u \right\|_{1,\Sigma}
\end{equation}
where we define
\[ \left\| \Tsig u \right\|^2_{1,\Sigma} := \left\| u \ResSig
\right\|^2_{H^1 (\Sigma ; g )} + \left\| \partial_t u
  \ResSig \right\|^2_{L^2 (\Sigma ; \chardsig )} \, .\]
\end{theorem}
Lars H\"ormander's proof can be extended with minor modifications
to the case where $g$ is in ${\cal C}^2 (\tx )$ and the coefficients of
$L_1$ are in $W^{1,\infty}_\mathrm{loc} (\tx )$~: this guarantees the
existence of ``regular'' solutions living in $H^2_\mathrm{loc} (\tx )$
which is enough for proving the energy estimates~; the whole proof can
then be reproduced using such solutions, instead of the ${\cal
  C}^\infty$ solutions used in the smooth case, to approach finite
energy solutions. This however is not quite enough for meeting the
standards imposed by H\"ormander in his final remark.

\section{Main results} \label{MainResults}

\subsection{The Cauchy problem}
\label{SectionCauchy}

We work on the same geometrical background but we now merely assume
the following~:
\begin{description}
\item[{\bf (H1)}] the metric $g$ is in ${\cal C}^0 (\tx ) \cap
W^{1,\infty}_\mathrm{loc} (\tx )$ and satisfies Property
\ref{ControlMetric}, the coefficients of $L_1$ are in
$L^\infty_\mathrm{loc} (\tx )$.
\end{description}
\begin{remark}
The regularity of $\dnu$ does not need to be lowered since the two
operators $\square$ corresponding to two choices of $\dnu$~: $\dnu_1 =
\gamma_1 \d x$, $\dnu_2 = \gamma_2 \d x$, such that $\gamma_2 -
\gamma_1 \in W^{1,\infty} (X)$, differ only by a first order operator
with bounded coefficients. Hence, the difference between two choices
of density is hidden in a black box~: the operator $L_1$. In fact, if
we study the natural covariant wave equation on $\tilde{X}$, this
black box already hides the difference between the simplified
d'Alembertian $\square$ defined in (\ref{Box}) and the covariant
d'Alembertian associated with the Lorentzian metric $\d t^2 -g$~:
\[ \partial^2_t - \Delta_g = \partial_t^2 - \frac{1}{| \det g |^\half}
\partial_\alpha \left( | \det g |^\half g^{\alpha \beta}
  \partial_\beta \right) \, . \]
\end{remark}
Because of the lack of regularity of the
coefficients of the equation, it is more natural to abandon part of
the continuity in time of the solutions. We give a first existence and
uniqueness result for solutions that are simply
$L^\infty_\mathrm{loc}$ in time with values in $H^1 (X)$. Strikingly
enough, it is then very easy to show that such solutions are in fact
continuous with values in $H^1 (X)$.
\begin{theorem} \label{Cauchy}
We introduce the space
\[ \tilde{\cal F} = L^{\infty}_\mathrm{loc} \left( \R_t \, ; H^1(X)
\right) \cap {\cal C}^1 \left( \R_t \, ; L^2(X) \right) \, .\]
Under the hypothesis {\bf (H1)}, for any $\left( u_0 , u_1 \right) \in
H^1 (X) \oplus L^2 (X)$, for any $s \in \R$, equation (\ref{WaveEq})
has a unique solution $u \in \tilde{\cal F}$ such that
\[ u_{|_{t=s}} = u_0 \, ,~ \partial_t u_{|_{t=s}} = u_1 \, . \]
Moreover any solution of (\ref{WaveEq}) in $\tilde{\cal F}$ belongs to
$\cal F$. Therefore, we still denote by $\cal E$ the space of
solutions of (\ref{WaveEq}) in $\tilde{\cal F}$. The elements of $\cal
E$ satisfy energy estimate (\ref{EnEst1}).
\end{theorem}
The next result states that when the operator $L_1$ is homogeneous of
the first order and exactly cancels the first order terms of the
d'Alembertian, we can get more regular solutions. This will in
particular be crucial for the Goursat problem.
\begin{corollary} \label{HomogCauchy}
For a metric $g$ in ${\cal C}^0 (\tx ) \cap
W^{1,\infty}_\mathrm{loc} (\tx )$ that satisfies Property
\ref{ControlMetric}, we consider the equation
\begin{equation} \label{HomogWaveEq}
\partial_t^2 u - g^{\alpha \beta} \partial_\alpha \partial_\beta u = 0
\end{equation}
corresponding to (\ref{WaveEq}) with
\[ L_1 = \gamma^{-1} \partial_\alpha \left( \gamma g^{\alpha \beta}
\right) \partial_\beta \, .\]
Note that $g$ and $L_1$ then satisfy hypothesis {\bf (H1)}.
As a consequence of theorem \ref{Cauchy}, for any $(u_0 , u_1) \in H^1
(X) \oplus L^2 (X)$, (\ref{HomogWaveEq}) admits in $\tilde{\cal F}$ a
unique solution $u$ such that
\[ u_{|_{t=s}} = u_0 \, ,~ \partial_t u_{|_{t=s}} = u_1 \, , \]
and we have in fact $u \in {\cal F}$.
Moreover, if $(u_0 , u_1) \in H^2 (X) \oplus H^1 (X)$, then the
solution $u$ satisfies
\[ u \in \bigcap_{l=0}^2 {\cal C}^l (\R_t \, ;~H^{2-l} (X)) \, .\]
\end{corollary}

\subsection{The Goursat problem}
\label{SectionGoursat}

We give an extension of theorem \ref{thmHorm} for a
metric $g$ that is merely continuously differentiable on $\tilde{X}$
in the case where the hypersurface $\Sigma$ is fully
characteristic. More precisely, we assume~:
\begin{description}
\item[{\bf (H2)}] the metric $g$ is in ${\cal C}^1
(\tilde{X})$, the coefficients of the first order terms of $L_1$ are
continuous on $\tilde{X}$ and the coefficients of the zero-order terms
of $L_1$ are in $L^{\infty}_\mathrm{loc} (\tilde{X} )$.
\end{description}
The hypersurface $\Sigma$ is still defined by (\ref{Sigma}) where
$\varphi~: X \rightarrow \R$ is a Lipschitz function, but it is now
required to be fully null, that is
\begin{equation} \label{NullSigma}
g^{\alpha \beta} (x , \varphi (x) ) \partial_\alpha \varphi (x)
\partial_\beta \varphi (x) = 1 ~~~ \mathrm{almost~everywhere~on~}X\, .
\end{equation}
Contrary to what one may think, this actually makes things slightly
easier since the measure $\chardsig$ vanishes everywhere on $\Sigma$
and therefore the trace of $\partial_t u$ on $\Sigma$ is no longer
relevant, only the more easily controlled trace of $u$ plays a part in
the characteristic Cauchy problem. This is what allows us to extend
the results of theorem \ref{thmHorm} to the case of a ${\cal C}^1$
metric. We have the following theorem~:
\begin{theorem} \label{Goursat}
Under the assumptions {\bf (H1)} , {\bf (H2)} and (\ref{NullSigma}) the
application
\[ \begin{array}{cccl} {\Tsig ~: }& {\cal E} & \longrightarrow &
  {H^1 (\Sigma  )} \\  & u & \longmapsto & {u\ResSig}
  \end{array} \]
is well defined and is an isomorphism.
\end{theorem}
The theorem is actually valid for slightly less regular metric and coefficients
of the first order terms. This becomes clear towards the end of the proof, in
the only part where we really need more than the minimum regularity setting
proposed by H\"ormander. The arguments are detailed in remark
\ref{RemProofCorol}.
\begin{theorem} \label{TheorMinReg}
The result of theorem \ref{Goursat} is still valid under the assumptions {\bf
(H1)}, (\ref{NullSigma}) and
\begin{description}
\item[{\bf ($\widehat{\mathbf{H2}}$)}] the metric $g$ is in
$L^\infty_\mathrm{loc} (\R_t \, ;~ {\cal C}^1 (X)) \cap
W^{1,\infty}_\mathrm{loc} (\R_t \, ;~ {\cal C}^0 (X)) $, the coefficients of the
first order terms of $L_1$ are in $L^\infty_\mathrm{loc} (\R_t \, ;~ {\cal C}^0
(X))$ and the coefficients of the zero-order terms of $L_1$ are in
$L^{\infty}_\mathrm{loc} (\tilde{X} )$.
\end{description}
\end{theorem}
The following additional result is implicit in the proof of theorem
\ref{Goursat}~:
\begin{corollary} \label{CorThmGoursat}
Let $u \in {\cal E}$, we define $v$ as
\[ \begin{array}{ccrcl} v & : & \R & \longrightarrow & {H^1 (\Sigma )
    \, , }  \\ & & t & \longmapsto & {u_{|_{\Sigma_t}} \, . }
\end{array} \]
Then $v \in {\cal C} \left( \R_t \, ; H^1 (\Sigma ) \right)$. This can
be expressed using the parametrization (\ref{Sigmat}) as follows~:
\[ v(t,x) := u(t+\varphi (x) , x) \, ,~v\in \left( \R_t \, ; H^1 (X )
\right) \, .\]
\end{corollary}

\section{Proofs of the main results}

\subsection{Proof of theorem \ref{Cauchy}}

We shall use the
following notations~: for $-\infty < t_1 < t_2 < +\infty$, ${\cal
  U}_{t_1 , t_2 } := ]t_1 , t_2 [ \, \times X$ and for $T>0$,
$\Omega_T := {\cal U}_{-T , T }$.

$\bullet$ {\bf First step~: uniqueness.} We show that the
energy estimate (\ref{EnEst1}) is valid for solutions of
(\ref{WaveEq}) in $\tilde{\cal F}$. We start by establishing an energy
estimate for all smooth functions, not assumed to satisfy
(\ref{WaveEq}), on $\tilde{X}$. For $v \in {\cal C}^\infty (\tx )$, if
we multiply $\square v + L_1 v$ by $\partial_t v$ and integrate on
$\Omega_T$ for $T>0$ given, the regularity of $v$
allows us to integrate by parts on $\Omega_T$ and to obtain
\begin{equation} \label{SmoothIntEst}
E (t,v ) \leq E(s,v) + 2 \int_{\Omega_T} \left| \partial_t v \right|
\, \left| \square v + L_1 v \right| \, \d t \, \dnu + K_1 (T,g,L_1 )
\int_{]s,t[} E(\tau , v) \d \tau
\end{equation}
where $K_1 (T,g,L_1 )$ is the continuous positive function of $T>0$,
the norms in $W^{1,\infty} (\Omega_T)$ of $g$ and $g^{-1}$ and the
norms of the coefficients of $L_1$ in $L^{\infty} (\Omega_T)$,
appearing in (\ref{EnEst1}). By density, estimate (\ref{SmoothIntEst})
carries over to functions $v$ in $H^2_\mathrm{loc} (\tx )$, but it is
not obvious that it remains valid on $\tilde{\cal F}$ because, using a
naive approximation, we cannot make sense either of the convergence of
the energy at times $s$ and $t$ or of the convergence of the term
containing the d'Alembertian. However, we show that any $u \in {\cal
  E}$ can be approached by a sequence $\{ u_k \}_k$ of more regular
functions such that estimate (\ref{SmoothIntEst}) for $u_k$ gives, as
$k \rightarrow +\infty$, estimate (\ref{EnEst1}) for $u$. There are
three constraints in the construction of the sequence
$\{ u_k \}_k$~:
\begin{enumerate}
\item we must have $E (t, u_k ) \rightarrow E(t,u)$ at least for
  almost every $t$~;
\item each $u_k$ has to belong to $H^2_\mathrm{loc} (\tx )$ so as to
  satisfy (\ref{SmoothIntEst})~;
\item $\square u_k + L_1 u_k$ must tend to zero at least weakly in
  $L^2_\mathrm{loc} (\tx )$ in order to obtain
\[ \int_{\Omega_T} \left| \partial_t u_k \right| \, \left| \square u_k
  + L_1 u_k \right| \, \d t \, \dnu \, \longrightarrow 0 \,
  ,~\mathrm{as}~k \rightarrow + \infty \, . \]
\end{enumerate}
The first constraint suggests to regularize $u$ in space only. The two
other constraints will then automatically be satisfied as well. In the
following proposition, we prove the existence of such a regularization
in trivial topology, i.e. on $\R_t \times \R^n$. Then, we use this
result locally on $\tx$ to construct the sequence $\{ u_k \}_k$.
\begin{proposition} \label{RegTrivTop}
We consider on $\R^n$ a time dependent metric $h(t)$ such that, for
all $t_1$, $t_2$, $-\infty < t_1 <t_2 < +\infty $,
\begin{eqnarray}
& h \in {\cal C}^0 \left( \R_t \times \R^n \right) \cap W^{1,\infty}
\left( ]t_1 , t_2 [ \, \times \R^n \right) \, , \label{Prop1h} \\
& \exists A,B \, ,~ 0 < A < B <+\infty \, ; ~A \id_n \leq h(t,x ) \leq
B \id_n \, , ~\forall (t,x) \in [t_1 , t_2 ] \, \times \R^n \,
, \label{Prop2h}
\end{eqnarray}
where $h$ is identified with its matrix in cartesian coordinates. We
introduce the operator
\[ \square_h = \partial^2_t - \partial_\alpha \left( h^{\alpha \beta }
  \partial_\beta \right) \, . \]
Let
\[ w \in L^\infty_\mathrm{loc} \left( \R_t \, ; H^1 (\R^n ) \right)
\cap {\cal C}^1 \left( \R_t \, ; L^2 (\R^n )\right) \]
such that
\[ \square_h w \in L^2 \left( ]t_1 , t_2 [ \, \times \R^n \right)
~\forall t_1,t_2 \, ;~ -\infty < t_1 <t_2 < +\infty \, .\]
If we consider on $\R^n$ a regularizing sequence (also called
approximate identity) defined in the usual manner
\[ \rho \in {\cal C}^\infty \left( \R^n \right) \, , ~ \supp \rho
\subset \bar{B} (0,1) \, ,~\rho \geq 0 \, ,~\int_{\R^n} \rho (x) \d x
= 1 \, ,~ \rho_k (x) := k^n \rho (kx) ~\forall k \in \N^* \, ,\] 
then the sequence $\{ w_k \}_k$ defined by convolution of $w$ with
$\rho_k$ over $\R^n$~:
\[w_k (t,x) = \left( w(t) * \rho_k \right) (x) \, ,\]
satisfies
\begin{description}
\item[(i)] $w_k (t) \rightarrow w (t)$ in $H^1 (\R^n )$ for all $t$
  such that $w(t) \in H^1 (\R^n )$~;
\item[(ii)] $w_k \rightarrow w$ in $L^p_\mathrm{loc} \left( \R_t \, ;
    H^1 (\R^n )\right) \cap {\cal C}^1 \left( \R_t \, ; L^2 (\R^n
    )\right)$ for all $1\leq p <+\infty$~;
\item[(iii)] the sequence $\left\{ \square_h w_k \right\}_k$ is bounded
  in $L^2_\mathrm{loc} \left( \R \times \R^n \right)$~;
\item[(iv)] $w_k \in H^2_\mathrm{loc} \left( \R \times \R^n \right)$
  $\forall k \in \N^*$.
\end{description}
\end{proposition}
{\bf Proof of proposition \ref{RegTrivTop}.}
\begin{description}
\item[(i) {\rm and} (ii)] are standard.
\item[(iii)] Since $\square_h w \in L^2_\mathrm{loc} \left( \R^{n+1}
  \right)$, we have
\[ \left( \square_h w \right) * \rho_k \rightarrow \square_h w ~
\mathrm{in}~ L^2_\mathrm{loc} \left( \R^{n+1} \right) \simeq
L^2_\mathrm{loc} \left( \R_t \, ; L^2_\mathrm{loc} \left( \R^n
  \right) \right) \, . \]
Hence, $ \left( \square_h w \right) * \rho_k$ is bounded in this
space. We consider
\begin{eqnarray}
\left( \square_h w \right) * \rho_k - \square_h w_k &=& \left(
  \partial_\alpha h^{\alpha \beta} \right) \left[ \left(
  \partial_\beta w \right) * \rho_k \right] + h^{\alpha \beta} \left[
  \left( \partial_\beta w \right) * \left( \partial_\alpha \rho_k
  \right) \right] \nonumber \\
&& - \left( h^{\alpha \beta} \partial_\beta w \right) *
  \left( \partial_\alpha \rho_k \right) \, . \nonumber
\end{eqnarray}
The first term is clearly bounded in $L^2_\mathrm{loc} \left( \R^{n+1}
\right)$ since
\[ \partial_\alpha h^{\alpha \beta} \in L^\infty \left( ]t_1 , t_2
  [\times \R^n \right) ~\mathrm{for~any} \, -\infty < t_1 < t_2 < +\infty
  \, ,~\partial_\beta w \in L^2_\mathrm{loc} \left( \R_t
  \, ; L^2\left( \R^n \right) \right)\]
and therefore
\[ \left( \partial_\alpha h^{\alpha \beta} \right) \left[ \left(
    \partial_\beta w \right) * \rho_k \right] \rightarrow \left(
    \partial_\alpha h^{\alpha \beta} \right) \partial_\beta w
    ~\mathrm{in}~L^2_\mathrm{loc} \left( \R_t \, ; L^2\left( \R^n
    \right) \right) \, .\]
Now
\begin{eqnarray}
&&\left\{ h^{\alpha \beta} \left[ \left( \partial_\beta w \right) *
  \left( \partial_\alpha \rho_k \right) \right] - \left( h^{\alpha
  \beta} \partial_\beta w \right) * \left( \partial_\alpha \rho_k
  \right) \right\} (t,x) \label{diff}  \\
&&= \int_{\supp \, \rho_k} \left\{ \left( h^{\alpha \beta} (t,x) -
  h^{\alpha \beta} (t,x-y) \right) \partial_\beta w (t,x-y) \,
  \partial_\alpha \rho_k (y) \right\} \d y \, . \nonumber
\end{eqnarray}
Denoting this quantity $F_k (t,x)$, we can estimate it as follows~:
for $-\infty < t_1 <t_2 < +\infty$, denoting $\Omega = ]t_1 , t_2 [
\times \R^n$, we have for $(t,x) \in \Omega$
\[ \left| F_k (t,x) \right| \leq \frac{1}{k} \left\| h^{-1}
  \right\|_{W^{1,\infty} (\Omega )} \sum_{\alpha , \beta} \int_{\supp
  \, \rho_k} \left| \partial_\beta w \left( t, x- y \right) \right|
  k^{n+1} \left| \partial_\alpha \rho (ky )\right| \d y \]
and putting $\chi_{\alpha , k} (y) = k^n \chi_\alpha (ky)$,
$\chi_\alpha (y) = \left| \partial_\alpha \rho (y) \right|$,
\[ \left| F_k (t,x) \right| \leq \left\| h^{-1} \right\|_{W^{1,\infty}
  (\Omega )} \sum_{\alpha , \beta} \left( \left| \partial_\beta w (t)
  \right| * \chi_{\alpha ,k} \right) (x) \, . \]
Using $\| \chi_{\alpha , k} \|_{L^1 (\R^n)} = \| \chi_{\alpha} \|_{L^1
  (\R^n)}$, we obtain
\[ \left\| F_k (t) \right\| _{L^2 (\R^n)} \leq \left\| h^{-1}
  \right\|_{W^{1,\infty} (\Omega )} \sum_{\alpha , \beta} \left\|
    \partial_\beta w (t) \right\|_{L^2 (\R^n )} \left\| \chi_{\alpha}
    \right\|_{L^1 (\R^n)} \]
and therefore
\[ \left\| F_k \right\| _{L^2 (\Omega )} \leq C \left\| h^{-1}
\right\|_{W^{1,\infty} (\Omega )} \left\| w \right\|_{H^1 (\Omega )}
\, , \]
where $C$ depends only on $\rho$. This proves {\it (iii)}.
\item[(iv)] We have $w \in {\cal C}^1 \left( \R_t \, ; L^2 (\R^n )
  \right)$, hence, for each $k\in \N^*$, $w_k$ is in ${\cal C}^1
  \left( \R_t \, ; {\cal C}^\infty (\R^n ) \right)$. Besides, we have
  proved that
\[ \square_h w_k = \partial^2_t w_k - \partial_\alpha \left( h^{\alpha \beta }
  \partial_\beta w_k \right) \in L^2_\mathrm{loc} \left( \R^{n+1}
  \right) \, .\]
Also
\[ w_k \in {\cal C}^1 \left( \R_t \, ; {\cal C}^\infty (\R^n )
\right) ~\mathrm{and}~ h^{\alpha \beta} \in W^{1,\infty } \left( ]t_1
  , t_2 [ \, \times \R^n \right) ~\forall ~-\infty < t_1 <t_2 <
  +\infty \]
entail
\[ \partial_\alpha \left( h^{\alpha \beta } \partial_\beta w_k \right)
\in L^2_\mathrm{loc} \left( \R_t \, ; L^2_\mathrm{loc} \left( \R^{n}
  \right) \right) \, .\]
Therefore
\[ \partial^2_t w_k \in L^2_\mathrm{loc} \left( \R^{n+1} \right) \]
which proves {\it (iv)} and concludes the proof of proposition
\ref{RegTrivTop}. \qed
\end{description}
We now proceed to constructing the sequence $\{ u_k \}_k$. We
consider~:
\begin{itemize}
\item $\left\{ \Omega^i \right\}_{1\leq i \leq N}$ a covering of $X$
  by open sets of trivial topology~;
\item $\left\{ {\cal U}^i \right\}_{1\leq i \leq N}$ a covering of $X$
  by open sets such that $\overline{{\cal U}^i} \subset \Omega^i$~;
\item $\left\{ \chi^i \right\}_{1\leq i \leq N}$ a partition of unity
  associated with $\left\{ {\cal U}^i \right\}$, i.e.
\[ \chi^i \in {\cal C}^\infty (X) \, ,~\supp \chi^i \subset
\overline{{\cal U}^i} \, ,~0\leq \chi^i \leq 1 \, ,~\sum_{i=1}^{N}
\chi^i = 1 ~\mathrm{on}~X \, . \]
\end{itemize}
Putting $v^i = \chi^i u$, we clearly have $v^i \in \tilde{\cal F}$ and
\[ \square v^i = \chi^i \square u - \gamma^{-1} \partial_\alpha \left(
  \gamma g^{\alpha \beta} \left( \partial_\beta \chi^i \right) u
  \right) - g^{\alpha \beta} \left( \partial_\alpha \chi^i \right)
  \left( \partial_\beta u \right) \, \in L^2_\mathrm{loc} \left(
  \tilde{X} \right) \, .\]
For $1 \leq i \leq N$, $\Omega^i$ is a trivial topology open subset of
$X$~; it can therefore be identified, by means of a global coordinate
system, with a bounded open set $\Omega^i$ in $\R^n$. The metric
$g_{|_{\R_t \times \Omega^i}}$ can be extended as a function $\ih
(t,x)$ on $\R_t \times \R^n$ satisfying (\ref{Prop1h}) and
(\ref{Prop2h}). The functions $v^i$ and $\square v^i$ on $\tx$ have
their support in $\R_t \times \overline{{\cal U}^i} \subset \R_t
\times \Omega^i$ and can therefore be considered as functions on $\R_t
\times \R^n$. Then, $v^i$ and $$\square_{\ih} v^i = \square v^i +
\mathit{~first~order~terms~in~}L^2_\mathrm{loc} \left( \R_t \, ; L^2
  \left( \Omega^i \right) \right) $$ satisfy the hypotheses of
proposition \ref{RegTrivTop}. Hence, for each $1\leq i\leq N$, we can
construct the sequence
\[ v^i_k = v^i * \rho_k \, , ~k \in \N^* \]
and it will satisfy properties {\it (i)}-{\it (iv)}. In addition,
\[ \supp v^i_k \subset \R_t \times \left( \overline{{\cal U}^i} +
\bar{B} \left( 0 , \frac{1}{k} \right) \right) ~\mathrm{in} ~ \R
\times \R^n \]
whence for $k$ large enough, $\square_{\ih} v^i_k$ and $v^i_k$ have
their support in $\R_t \times \Omega^i$ and can be
considered as functions on $\tx$. In this manner, for a given $\bar{k}
\in \N^*$, we obtain a sequence $\{ u_k \}_{k \geq \bar{k}}$ defined
by
\[ u_k = \sum_{i=1}^{N} v^i_k \]
such that
\begin{eqnarray}
& u_k \in H^2_\mathrm{loc} (\tx )\, ,\label{Reguk} \\
& u_k (t) \rightarrow u(t) ~\mathrm{in}~H^1 (X) \, ,~\forall
t~\mathrm{such~that~} u(t)\in H^1(X) \, , \label{Convuk1} \\
& u_k \rightarrow u ~\mathrm{in}~ L^p_\mathrm{loc} \left( \R_t \, ;
  H^1 (X) \right) \cup {\cal C}^1 \left( \R_t \, ; L^2 (X) \right)
~\forall p\, ,\, 1\leq p <+\infty \, , \label{Convuk2} \\
& \square u_k + L_1 u_k ~\mathrm{bounded~in~} L^2_\mathrm{loc} (\tx )
\, . \label{Bounduk}
\end{eqnarray}
Note that (\ref{Bounduk}) is an easy consequence of property {\it
  (iii)} for $v_i^k$ since
\[ \square u_k + L_1 u_k - \sum_{i=1}^{N} \square_{\ih} v^i_k \]
is a sum of first or zero order derivatives of the $v^i_k$ with
coefficients in $L^\infty_\mathrm{loc} (\tx )$ and these terms
converge in $L^2_\mathrm{loc} (\tx )$ by (\ref{Convuk2}). By
(\ref{Reguk}), each $u_k$ satisfies (\ref{SmoothIntEst})~: for all
$T>0$, for all $t,s\in [-T ,T]$,
\begin{equation} \label{Estuk}
E(t,u_k) \leq E(s,u_k) + 2 \int_{\Omega_T} \left| \partial_t u_k
\right| \left| \square u_k + L_1 u_k \right| \d t \d \nu + K_1
(T,g,L_1 ) \int_{]s,t[} E(\tau , u_k) \d \tau \, .
\end{equation}
Properties (\ref{Convuk1}) and (\ref{Convuk2}) imply that for almost
all $s,t \in [-T,T]$ (more precisely for all $s,t$ such that both
$u(s)$ and $u(t)$ belong to $H^1 (X)$)
\[ E(t,u_k ) \rightarrow E(t,u ) ~\mathrm{and}~E(s,u_k ) \rightarrow
E(s,u ) \, .\]
Property (\ref{Convuk2}) also entails the convergence of the last term
of the inequality
\[ \int_{]s,t[} E(\tau ,u_k ) \d \tau \rightarrow \int_{]s,t[} E(\tau
,u ) \d \tau \, .\]
Extracting a subsequence if necessary, (\ref{Bounduk}) entails that
$\square u_k +L_1 u_k$ converges weakly in $L^2 (\Omega_T )$, the
limit being zero since, by (\ref{Convuk1}) and using the fact that $u$
is a solution of (\ref{WaveEq}), $\square u_k +L_1 u_k$ converges
towards zero in ${\cal D}' (\Omega_T )$, the space of distributions on
$\Omega_T$. Besides, $\partial_t u_k$
converges strongly towards $\partial_t u$ in $L^2 (\Omega_T )$, whence
\[ \int_{\Omega_T} \left| \partial_t u_k \right| \left| \square u_k +
  L_1 u_k \right| \d t \d \nu \longrightarrow 0 \, .\]
Consequently, for all $t,s\in [-T,T]$ such that $u(t),u(s)\in H^1
(X)$,
\begin{equation} \label{EnEstFtilde}
E(t,u) \leq E(s,u ) + K_1 (T,g,L_1 ) \int_{]s,t[} E(\tau , u ) \d
\tau \, .
\end{equation}
This gives (\ref{EnEst1}) for $u$ in the following sense~: if $u(s)
\in H^1 (X )$, $s \in [-T ,T]$, then for almost all $t \in [-T ,T]$,
\[ E(t,u ) \leq E(s,u ) e^{K_1 (T,g,L_1 ) |t-s|} .\]
The uniqueness of solutions to the Cauchy problem for (\ref{WaveEq})
in $\tilde{\cal F}$ follows. \qed

$\bullet$ {\bf Second step~: existence.} Let $( \phi , \psi )
\in H^1 (X) \oplus L^2 (X)$, $s\in \R$, we wish to find $u \in
\tilde{\cal F}$ such that
\begin{equation} \label{CauchyProb}
\square u +L_1 u =0 \, ,~u_{|_{t=s}} = \phi \, , ~\partial_t
u_{|_{t=s}} = \psi \, .
\end{equation}
In order to use the well posedness of the Cauchy problem in the smooth
case, we regularize the metric and the coefficients of $L_1$. For $k
\in \N^*$, we define
\begin{description}
\item[(a)] $\gk$ a time-dependent riemannian metric on $X$, $\gk
  \in {\cal C}^\infty (\tx )$~;
\item[(b)] $\oplk = \bk^0 \partial_t + \bk^\alpha \partial_\alpha
  + \ck \, ,~\bk^0\, ,~\bk^\alpha \, ,~\ck \, \in {\cal C}^\infty (\tx )$~;
\end{description}
such that~:
\begin{eqnarray}
& \gk \longrightarrow g ~\mathrm{in}~{\cal C}^0 (\tx) \cap
H^{1}_\mathrm{loc} (\tx ) \, ; \label{Convgk} \\
& \gk ~\mathrm{bounded~in~} W^{1,\infty}_\mathrm{loc} (\tx ) \,
; \label{Boundgk} \\
& \left. \begin{array}{ccc} \ck & \longrightarrow & c \\ \bk^0 &
  \longrightarrow & b^0 \\ \bk^\alpha & \longrightarrow & b^\alpha
\end{array} \right\} ~\mathrm{in}~L^p_\mathrm{loc} (\tx )~\forall p\,
;~1\leq p < +\infty \, ; \label{ConvL1k} \\
& \bk^0  , \bk^\alpha  ,\ck
~\mathrm{bounded~in~} L^\infty_\mathrm{loc} (\tx ) \label{BoundL1k}
\end{eqnarray}
and there exist two positive continuous functions $D_1$ and $D_2$ such
that
\begin{equation} \label{UpLowBoundgk}
\forall k \in \N^* \, ,~\forall (t,x) \in \tx \, ,~D_1 (t) \id_n \leq
\gk_{\alpha \beta} (t,x) \leq D_2 (t) \id_n \, .
\end{equation}
Typically, such sequences are constructed using coordinate charts and,
in each domain, convolution by a regularizing sequence on $\R_t \times
\R^n$~; this is similar to what we did for constructing the sequence
$u_k$ in the first step of the proof, but now, the regularizing
sequence and the convolution involve time as well as space variables.

For each $k$, we consider the equation
\begin{equation} \label{RegulEq}
\frac{\partial^2 v}{\partial t^2}- \gamma^{-1} \frac{\partial}{\partial
  x^\alpha} \left( \gamma \left( \gk^{\alpha \beta}\right)
  \frac{\partial v}{\partial x^\beta} \right) + \oplk v = 0 \, .
\end{equation}
Theorem \ref{thmHorm} tells us that (\ref{RegulEq}) has a unique
solution $v_k \in {\cal F}$ such that $v_k(s) = \phi$ and $\partial_t
v_k (s) = \psi$. This solution satisfies the energy estimate
\begin{equation} \label{EnEstvk}
\forall T > |s| \, ,~\forall t \in [-T,T] \, ,~E_k (t,v_k ) \leq E_k
(s,v_k ) \, e^{K_1 (T, \gk , \oplk ) |t-s|} ,
\end{equation}
where $E_k$ is the energy (\ref{Energy}) defined using the metric
$\gk$ instead of $g$. $E_k (t,.)$, just like $E(t,.)$, is (uniformly
in $k$ and locally uniformly in time) equivalent to the norm in $H^1
(X) \oplus L^2 (X)$. Besides, (\ref{Boundgk}) and (\ref{BoundL1k})
imply that $\left\{ K_1 (T, \gk , \oplk ) \right\}_{k}$ is bounded in
$\R^+$. The upshot of all this is that $\{ v_k \}_k$ is bounded in
$\cal F$. Hence, for $T>|s|$ fixed, extracting a subsequence if
necessary, we can assume the convergence of $\{ v_k \}$
in the following spaces (we call $u$ the common limit)~:
\begin{eqnarray}
& v_k \rightarrow u ~\mathrm{in}~H^1 (\Omega_T )-w \, ,
\label{ConvH1vk} \\
& v_k \rightarrow u ~\mathrm{in}~H^\mu (\Omega_T ) \, ,~\forall \mu <1
\, , \label{ConvHmuvk}
\end{eqnarray}
where ``$-w$'' denotes the weak topology. Hence, by standard trace
theorems
\begin{eqnarray}
& v_k \rightarrow u ~\mathrm{in}~ {\cal C} \left( [-T,T] \, ; L^2 (X)
\right) \label{ConvTracevk}
\end{eqnarray}
and by the Banach-Alaoglou theorem
\begin{eqnarray}
& v_k \rightarrow u ~\mathrm{in}~ L^\infty \left( ]-T,T[ \, ; H^1 (X)
\right) -w-* \, , \label{ConvBA1vk} \\
& \partial_t v_k \rightarrow \partial_t u ~\mathrm{in}~ L^\infty
\left( ]-T,T[ \, ; L^2 (X) \right) -w-* \, , \label{ConvBA2vk}
\end{eqnarray}
where ``$-w-*$'' denotes the weak star topology. Now the convergences
(\ref{ConvH1vk}), (\ref{Convgk}) and (\ref{ConvL1k}) imply
\[ \left. \begin{array}{rcl} {\partial^2_t v_k} & \longrightarrow &
    {\partial^2_t u} \\ \\{\gamma^{-1} \partial_\alpha \left( \gamma
    \left( \gk^{\alpha \beta } \right) \partial_\beta v_k \right)} &
    \longrightarrow & {\gamma^{-1} \partial_\alpha \left( \gamma
    g^{\alpha \beta } \partial_\beta u \right)} \\ \\
    {\oplk v_k} & \longrightarrow & {L_1 u} \end{array} \right\}
    ~\mathrm{in}~{\cal D}' (\Omega_T ) \, ,\]
whereby $u$ satisfies equation (\ref{WaveEq}) in the sense of distributions on
$\Omega_T$. Using uniqueness, we have thus constructed a solution $u$ of
(\ref{WaveEq}) defined on $\tx$ and that belongs to $\tilde{\cal F}$. Indeed, we
know that
\begin{equation} \label{Regulu1}
u \in L^\infty_\mathrm{loc} \left( \R_t \, ; H^1 (X) \right) \, ,~
\partial_t u \in L^\infty_\mathrm{loc} \left( \R_t \, ; L^2 (X)
\right) \, .
\end{equation}
Since $g \in W^{1,\infty}_\mathrm{loc} (\tx )$ and the coefficients of
$L_1$ are in $L^\infty_\mathrm{loc} (\tx)$, (\ref{WaveEq}) entails
\begin{equation} \label{Regulu2}
\partial^2_t u = \gamma^{-1} \partial_\alpha \left( \gamma g^{\alpha
    \beta } \partial_\beta u \right) - L_1 u \in L^\infty_\mathrm{loc}
    \left( \R_t \, ; H^{-1} (X) \right) \, .
\end{equation}
Using J.-L. Lions's principle of intermediate derivatives,
(\ref{Regulu1}) and (\ref{Regulu2}) imply
\[ \partial_t u \in {\cal C} \left( \R_t \, ; L^2 (X) \right) \]
and therefore $u \in \tilde{\cal F}$. The last things to check are the
two initial data conditions. The initial value of $u$ is easy~; using
(\ref{ConvTracevk})
\[ v_k (s) = \phi \longrightarrow u(s) ~\mathrm{in} ~L^2 (X) \]
whence $u(s) = \phi$. The trace of $\partial_t u$ at $t=s$ requires
more care. We write
\[ \partial^2_t \left( u -v_k \right) = \gamma^{-1} \partial_\alpha
\left( \gamma \left[ \gk^{\alpha \beta } \partial_\beta v_k -
    g^{\alpha \beta } \partial_\beta u \right] \right) + \oplk v_k -
    L_1 u \, . \]
(\ref{Boundgk}), (\ref{BoundL1k}) and (\ref{ConvH1vk}) imply the
boundedness in $L^2 (]-T,T[\, ; H^{-1} (X) )$ of $\partial^2_t \left(
  u-v_k \right)$. This allows us, first, to write for $t\in [-T,T]$,
$t_0$ fixed in $[-T,T]$,
\begin{equation} \label{dtuContinu}
\partial_t u(t) - \partial_t v_k (t) = \partial_t u(t_0) -
\partial_t v_k (t_0) + \int_{]t_0 ,t[} \partial^2_t \left( u - v_k
\right) (\tau ) \d \tau \, ,
\end{equation}
second, extracting another subsequence if necessary, to assume
\[ \partial^2_t \left( u-v_k \right) \longrightarrow 0
~\mathrm{in}~L^2 \left( ]-T,T[\, ; H^{-1} (X) \right) -w \, . \]
This last convergence gives
\[ \int_{]t_0 ,t[} \partial^2_t \left( u - v_k \right) (\tau ) \d \tau
\longrightarrow 0~\mathrm{in}~L^2 \left( ]-T,T[\, ; H^{-1} (X) \right)
-w \, . \]
Since (\ref{ConvH1vk}) implies
\[ \partial_t v_k \rightarrow \partial_t u ~\mathrm{in}~L^2 (\Omega_T
) -w \hookrightarrow L^2 (]-T , T[ \, ;~H^{-1} (X) ) -w \, ,\]
we deduce from (\ref{dtuContinu}) that
\[ \partial_t v_k (t_0 ) \longrightarrow \partial_t u(t_0 )
~\mathrm{in}~H^{-1} (X) -w \, ,~\forall t_0 \in [-T,T] \, . \]
In particular, for $t_0 =s$,
\[ \Psi = \partial_t v_k (s) \longrightarrow \partial_t u(s)
~\mathrm{in}~H^{-1} (X)-w \]
which gives us $\partial_t u(s) =\Psi$ and concludes the second part
of the proof of theorem \ref{Cauchy}. \qed

$\bullet$ {\bf Third step~: continuity in time of the
solutions.} We consider $u$ the unique solution  in $\tilde{\cal F}$
of the Cauchy problem (\ref{CauchyProb}), $T>|s|$ and $v_k$ the
sequence constructed in the second step of the proof. For any fixed
$t$ in $[-T,T]$, the energy estimate (\ref{EnEstvk}) implies that $\{
v_k (t)\}_k$ is bounded in $H^1 (X)$. Hence, extracting a subsequence
if necessary, we can assume that $v_k(t)$ converges
weakly in $H^1 (X)$. This together with the strong convergence
(\ref{ConvTracevk}) guarantees that $u(t)$ belongs to $H^1 (X)$. The
construction of the sequence $v_k$ can be made for any fixed
$T>|s|$. It therefore turns out that
\[ u(t) \in H^1 (X ) ~\forall t \in \R \]
and hence, the energy estimate (\ref{EnEstFtilde}) is valid for
all $t,s$. This implies in particular that for any solution $u$ of
(\ref{WaveEq}) in $\tilde{\cal F}$, the energy $E(t,u)$ is continuous
in time. Besides, it is easy to show that
\begin{equation} \label{WeakContFtilde}
u \in \tilde{\cal F} \Longrightarrow u \in {\cal C} \left( \R_t \,
  ; H^1 (X)-w \right) \, .
\end{equation}
The continuity of the energy therefore entails the strong continuity
of $u$ in time with values in $H^1 (X)$, which proves $u \in \cal
F$. We now prove (\ref{WeakContFtilde}). Let $v \in \tilde{\cal F}$
and $w \in H^1 (X)$. We put for $t\in \R$
\[ f(t) = \left< v(t) , w \right>_{H^1(X)} \, . \]
Given $t_0 \in \R$ we show the continuity of $f$ at $t_0$. Let $t_n
\rightarrow t_0$ and $\{ w_k \}_k$ a sequence in ${\cal C}^\infty (X)$
converging towards $w$ in $H^1 (X)$. For each $k$, using $v \in {\cal
  C} (\R_t \, ;~L^2(X))$, we have
\[  \left< v(t_n) -v(t_0 ) , w_k \right>_{H^1(X)} =  \left< v(t_n ) -
  v(t_0 ) , \left (1-\Delta_h \right) w_k \right>_{L^2 (X)}
  \longrightarrow 0 \, ,~n\rightarrow +\infty \, , \]
where $\Delta_h = \gamma^{-1} \partial_\alpha \left( \gamma h^{\alpha \beta }
\partial_\beta \right)$ is the Laplacian associated with the metric $h$ on $X$,
introduced in section \ref{SectionHormRes} to define the $H^1$ norm on $X$, and
to which the measure $\d \nu$ is associated. We write
\[ \left< v(t_n) -v(t_0 ) , w \right>_{H^1(X)} = \left< v(t_n) -v(t_0
    ) , w- w_k \right>_{H^1(X)} + \left< v(t_n) -v(t_0 ) , w_k
    \right>_{H^1(X)} \, .\]
Consider $\varepsilon >0$. Using the fact that $v \in
L^\infty_\mathrm{loc} (\R_t \, ;H^1 (X) )$, we choose $k$ large enough
so that for all $n$
\[ \left| \left< v(t_n) -v(t_0 ) , w- w_k \right>_{H^1(X)} \right|
\leq \varepsilon /2 \, ,\]
then, for this value of $k$, we choose $n$ large enough so that
\[ \left| \left< v(t_n) -v(t_0 ) , w_k \right>_{H^1(X)} \right| \leq
  \varepsilon /2 \, .\]
This proves the continuity of $f(t)$ and concludes the proof of
theorem \ref{Cauchy}. \qed

\subsection{Proof of corollary \ref{HomogCauchy}}

First, we write
(\ref{HomogWaveEq}) as a special case of (\ref{WaveEq})~:
\[ \partial_t^2 u - g^{\alpha \beta} \partial_\alpha \partial_\beta u
= \square u + \gamma^{-1} \partial_\alpha \left( \gamma g^{\alpha
    \beta} \right) \partial_\beta u =0\, .\]
Since $g, g^{-1} \in W^{1,\infty}_\mathrm{loc} (\tilde{X})$, the
coefficients of the first order operator clearly belong to
$L^\infty_\mathrm{loc} (\tilde{X})$ and we are in the framework of
theorem \ref{Cauchy}. To check that we can get more regular solutions,
we simply apply a partial derivation to (\ref{HomogWaveEq})~:
\[ \partial_\mu \left( \partial_t^2 u - g^{\alpha \beta}
  \partial_\alpha \partial_\beta u \right) = \left( \partial_t^2 -
  g^{\alpha \beta} \partial_\alpha \partial_\beta \right) \partial_\mu
u - \left( \partial_\mu g^{\alpha \beta} \right) \partial_\alpha
\partial_\beta u \, .\]
We can therefore write the following system of equations~:
\[ \left\{ \begin{array}{rcl}
{\left( \partial_t^2 - g^{\alpha \beta} \partial_\alpha \partial_\beta
  \right) u} & = & 0 \\
{\left( \partial_t^2 - g^{\alpha \beta} \partial_\alpha \partial_\beta
\right) \partial_\mu u} & = & {\left( \partial_\mu g^{\alpha
      \beta} \right) \partial_\alpha \partial_\beta u \, ,~\mu
=1,...,n.} \end{array} \right. \]
This system is of the form
\begin{equation} \label{DerivedSyst}
\left( \partial_t^2 - g^{\alpha \beta} \partial_\alpha \partial_\beta
  \right) U = L_1 U \, ,~U= \,^t\left( u, \partial_1 u , ... ,
    \partial_n u \right) \, ,
\end{equation}
where $L_1$ is a first order differential operator whose coefficients
belong to $L^\infty_\mathrm{loc} (\tilde{X})$ since the metric $g$ is
in $W^{1,\infty}_\mathrm{loc} (\tilde{X})$. By theorem \ref{Cauchy}
(in the case where the unknown function is a vector field), the system
(\ref{DerivedSyst}) admits a well-posed Cauchy problem in $\tilde{\cal
F}$ and the solutions belong to $\cal F$. This guarantees the
additional regularity of solutions of (\ref{HomogWaveEq}) for data in
$H^2 \oplus H^1$ and concludes the proof of corollary
\ref{HomogCauchy}. \qed

\subsection{Proof of theorem \ref{Goursat}.}

$\bullet$ {\bf Inequalities (\ref{EnEst2}) and (\ref{EnEst3}).} To prove
these inequalities with our regularity assumptions, we define a
regularization of the solution $u$ by functions $u_k$ in
$H^2_\mathrm{loc} (\tilde{X})$ that satisfy estimates of type
(\ref{EnEst2})-(\ref{EnEst3}), with constants uniform in $k$, and that
converge towards $u$ strongly in $H^1 (X_t)$ for all $t$ and in $H^1
(\Sigma )$. This makes a crucial use of corollary \ref{HomogCauchy}.

We write equation (\ref{WaveEq}) as follows
\begin{equation} \label{WEPerturbHomog}
\partial_t^2 u - g^{\alpha \beta} \partial_\alpha \partial_\beta u
+\tilde{L}_1 u = 0 \, ,
\end{equation}
where
\begin{gather*}
\tilde{L}_1 = - \gamma^{-1} \left( \partial_\alpha \left(
\gamma g^{\alpha \beta} \right) \right) \partial_\beta +L_1 =:
p^0\partial_t + p^\beta \partial_\beta + q \, ,\\
p^0 = b^0 \, ,~p^\beta = - \gamma^{-1} \left( \partial_\alpha \left(
\gamma g^{\alpha \beta} \right) \right) +b^\beta \, ,~q=c \, .
\end{gather*}
The coefficients of $\tilde{L}_1$ satisfy
\[ p^0\, ,~p^\beta \in {\cal C}^0 (\tilde{X}) \, ,~q \in
L^\infty_\mathrm{loc} ( \tilde{X}) \, .\]
We define an approximation of equation (\ref{WEPerturbHomog}) in which
only the coefficients of $\tilde{L}_1$ are regularized~:
\begin{equation} \label{ApproxL1k}
\partial_t^2 u - g^{\alpha \beta} \partial_\alpha \partial_\beta u
+\tilde{L}_1^k u = 0 \, , \tilde{L}_1^k = p^0_k \partial_t + p^\beta_k
\partial_\beta + q_k \, ,
\end{equation}
where the coefficients of $\tilde{L}_1^k$ satisfy
\begin{gather}
p_k^0 \, ,~p_k^\beta \, ,~q_k \in {\cal C}^\infty (\tilde{X}) \, ,
\label{RegCoeffs} \\
p_k^0 \rightarrow p^0~\mathrm{in}~{\cal C} (\tilde{X})\, ,~p_k^\beta
\rightarrow p^\beta~\mathrm{in}~{\cal C} (\tilde{X})\, ,
\label{ConvCoeffs1} \\
q_k \rightarrow q~\mathrm{in}~L^p_\mathrm{loc} (\tilde{X})\, ,~\forall
1\leq p <+\infty \, , \label{ConvCoeffs2} \\
q_k~\mathrm{bounded~in}~L^\infty (]-T , T [ \times X)~\forall T>0 \,
. \label{qboundLinfty}
\end{gather}
Let $(u_0 , u_1) \in H^1 (X) \oplus L^2 (X)$. We consider $u \in {\cal
  F}$ the solution of (\ref{WEPerturbHomog}) such that $u(0) = u_0$
and $\partial_t u (0) = u_1$. We also consider some sequences $\{
u_0^k \}_k$ and $\{ u_1^k \}_k$ of smooth functions on $X$ such that
\begin{equation} \label{ConvInData}
u_0^k \rightarrow u_0 ~\mathrm{in}~H^1(X) ~\mathrm{and} ~u_1^k
\rightarrow u_1 ~\mathrm{in}~L^2 (X)\, .
\end{equation}
Let $u_k \in H^2_\mathrm{loc} (\tilde{X})$ the solution in
$\cal F$ of (\ref{ApproxL1k}) such that $u_k (0) = u_0^k$ and
$\partial_t u_k (0) = u_1^k$. For each $k$, $u_k$ satisfies estimates of type
(\ref{EnEst2}) and (\ref{EnEst3}) uniformly in $k$, more precisely there exists
$C>0$ such that, for all $k$~:
\begin{eqnarray}
\left\| {u_k}_{|_\Sigma} \right\|^2_{1,\Sigma} \leq C \left( \left\| u_0^k
\right\|^2_{H^1 (X)} + \left\| u_1^k \right\|^2_{L^2(X)} \right) \label{Est1u_k}
\\
\left\| u_0^k \right\|^2_{H^1 (X)} + \left\| u_1^k \right\|^2_{L^2(X)} \leq C
\left\| {u_k}_{|_\Sigma} \right\|^2_{1,\Sigma} \label{Est2u_k}
\end{eqnarray}
In order to establish (\ref{EnEst2}) and (\ref{EnEst3}) for $u$, we only need to
prove that $u_k$ converges towards $u$ in $H^1 (\Sigma )$, since
\[\left\| u_0^k \right\|^2_{H^1 (X)} + \left\| u_1^k \right\|^2_{L^2(X)}
\rightarrow \left\| u_0 \right\|^2_{H^1 (X)} + \left\| u_1 \right\|^2_{L^2(X)} =
E(0,u) ~\mathrm{as~} k\rightarrow +\infty \, .\]
We will use the following proposition.
\begin{proposition} \label{3EstimatesRegSol}
We consider the equation
\begin{equation} \label{HomogSmoothPertModel}
\partial_t^2 v - g^{\alpha \beta} \partial^2_{\alpha \beta}v + L v = f
\end{equation}
where $L$ is a first order differential operator with smooth coefficients on
$\tilde{X}$ and the source $f$ belongs to $L^2_\mathrm{loc} (\tilde{X})$.
\begin{enumerate}
\item We consider $T> 0$ and $\Omega_T = ]-T,T[\times X$, there exists a
continuous positive function $C_1 (T, g,L)$ of $T$, the norms of $g$ and
$g^{-1}$ in $W^{1,\infty} (\Omega_T)$ and the norms of the coefficients of $L$
in $L^\infty (\Omega_T)$, such that, for any solution $v$ of
(\ref{HomogSmoothPertModel}) in $\cal F$ and for all $t,s \in [-T,T]$~:
\begin{equation}
E (t , v) \leq C_1 (T, g,L) \left( E(s , v) + \left\| f \right\|^2_{L^1
(]-T,T[\, ;~L^2 (X))} \right) \, . \label{Ineq0H2}
\end{equation}
\item We now consider $T>\max \{ |\min \varphi | \, ,~ | \max \varphi |\}$.
There exist continuous positive functions $C_2 (T,g,L)$, $C_3 (T,g,L,f)$ of $T$,
the norms of $g$ and $g^{-1}$ in $W^{1,\infty} (\Omega_T)$ and the norms of the
coefficients of $L$ in $L^\infty (\Omega_T)$, such that, for any solution $v$ of
(\ref{HomogSmoothPertModel}) in $H^2_\mathrm{loc} (\tilde{X})$ and for any $s\in
[-T,T]$, we have\footnote{The existence of such solutions is not guaranteed in
the general case because of the low regularity of $f$, but we will use this
proposition in cases where we know such solutions, namely the functions $u_k$ or
rather the difference $u_k - u_l$ between two such solutions~; see equation
(\ref{Equk-ul}).}~:
\begin{gather}
\left\| v_{|_\Sigma} \right\|_{H^1 (\Sigma )} \leq C_2 (T,g,L) \left( E(s,v) +
\left\| f \right\|^2_{L^1 (]-T,T[\, ;~L^2 (X))} \right) \, , \label{Ineq1H2}\\
E(s,v) \leq C_3 (T,g,L) \left( \left\|v_{|_\Sigma} \right\|_{H^1 (\Sigma )} +
\left\| f \right\|^2_{L^2 (\Omega_T)} \right) \, . \label{Ineq2H2}
\end{gather}
\end{enumerate}
\end{proposition}
\begin{remark}
Estimate (\ref{Ineq2H2}) will not be useful to us, we have given it for
completeness.
\end{remark}
{\bf Proof.}
\begin{enumerate}
\item {\bf Proof of (\ref{Ineq0H2}).}
We have obtained in the proof of theorem \ref{Cauchy} that estimate
(\ref{EnEst1}) is valid for solutions of (\ref{WaveEq}) in ${\cal F}$
under the assumption {\bf (H1)}. If we consider some source $f$ in
$L^1_\mathrm{loc} (\R_t \, ;~L^2(X))$, we still have existence and
uniqueness in $\cal F$ of the solutions of
\[ \square v +L_1 v= f\]
and these solutions are given in terms of their initial data at time
$s$ by the Duhamel formula\footnote{This is established by a standard
  fixed point argument.}
\begin{equation} \label{Duhamel}
\left( \begin{array}{c} {v(t)} \\ {\partial_t v(t)} \end{array}
\right) = {\cal U} (t,s) \left( \begin{array}{c} {v(s)} \\ {\partial_t
      v(s)} \end{array} \right) + \int_s^t {\cal U} (t,\tau ) \left(
  \begin{array}{c} {0} \\ {f(\tau )} \end{array} \right) \d \tau \, ,
\end{equation}
where ${\cal U} (t,s)$ denotes the propagator for equation
(\ref{WaveEq}), that to initial data $^t (u(s) , \partial_t u(s))$,
associates the solution at time $t$~: $^t (u(t) , \partial_t
u(t))$.
Equation (\ref{HomogSmoothPertModel}) in the source-free case can be
written as
\[ \square v + \gamma^{-1} \left( \partial_\alpha \left(
\gamma g^{\alpha \beta} \right) \right) \partial_\beta v + L v = 0 \,
\]
and therefore the solutions satisfy estimate (\ref{EnEst1}) with a bound $e^{K_1
(T,g,L)|t-s|}$ where $K_1$ is a continuous positive function of $T$, the norms
of $g$ and $g^{-1}$ in $W^{1,\infty} (]-T,T[\times X)$ and the norms of the
coefficients of $L$ in $L^\infty (]-T,T[\times X)$. This together with
(\ref{Duhamel}) entail (\ref{Ineq0H2}) with $C_1 = e^{2TK_1}$ for solutions in
$\cal F$ of equation (\ref{HomogSmoothPertModel}) with a source $f\in
L^1_\mathrm{loc} (\R_t \, ;~L^2(X))$.
\item {\bf Proof of (\ref{Ineq1H2}).}
The fact that we are dealing with a solution that is
locally $H^2$ allows us to use the same type of integrations by parts
as Lars H\"ormander. For $v \in H^2_\mathrm{loc} (\tilde{X})$ solution
of (\ref{HomogSmoothPertModel}), we write
\begin{eqnarray}
0 &=& 2 \partial_t v \left( \partial_t^2 v - g^{\alpha \beta}
  \partial^2_{\alpha \beta} v + L v -f\right) \nonumber \\
&=& \partial_t \left[ \left( \partial_t v \right)^2 + g^{\alpha \beta}
  \partial_\alpha v \partial_\beta v + v^2 \right] -2 \gamma^{-1}
\partial_\alpha \left[ \gamma g^{\alpha \beta} \partial_t v
    \partial_\beta v \right] \nonumber \\
&& + 2 \partial_t v L v + 2 \gamma^{-1} \partial_\alpha \left( \gamma
  g^{\alpha \beta} \right) \partial_t v \partial_\beta v - \left(
  \partial_t g^{\alpha \beta} \right) \partial_\alpha v \partial_\beta
v -2 v \partial_t v - 2 f\partial_t v \label{EnergyDensity}
\end{eqnarray}
Integrating (\ref{EnergyDensity}) on the domain $\Omega_T^- =\{ -T
\leq t \leq \varphi (x) \}$ for the measure $\d t \d \nu = \gamma \d t
\d x$, we obtain
\begin{eqnarray*}
0 &=& \int_\Sigma \left( \left( \partial_t v \right)^2 + g^{\alpha
    \beta} \partial_\alpha v \partial_\beta v + v^2 \right) \dsig -
\int_{X_{_{-T}}} \left( \left( \partial_t v \right)^2 +
  g^{\alpha \beta} \partial_\alpha v \partial_\beta v + v^2 \right) \d
\nu \\
&&+ \int_\Sigma 2 g^{\alpha \beta} \partial_t v \partial_\alpha
\varphi \partial_\beta v \dsig \\
&&+ \int_{\Omega_T^-} \left( 2
  \partial_t v L v + 2 \gamma^{-1} \partial_\alpha \left( \gamma
    g^{\alpha \beta} \right) \partial_t v \partial_\beta v \right. \\
&& \hspace{0.5in} \left. - \left( \partial_t g^{\alpha \beta} \right)
  \partial_\alpha v \partial_\beta v -2 v \partial_t v - 2 f\partial_t
  v \right) \d t\d \nu \, .
\end{eqnarray*}
The first three terms give (using the fact that $\Sigma$ is totally
null)
\begin{gather*}
-E(-T ,v) + \int_\Sigma \left( \partial_t v\right)^2 \chardsig +
\int_\Sigma \left( g^{\alpha \beta} \left(
    \partial_\alpha v + \partial_\alpha \varphi \partial_t v \right)
  \left( \partial_\beta v + \partial_\beta \varphi \partial_t v
  \right) + v^2 \right) \dsig \\
= -E(-T,v) + \left\| v_{|_\Sigma} \right\|^2_{H^1 (\Sigma )}
\end{gather*}
and the other terms, thanks to the assumptions on $g$ and the
coefficients of $L$, can be estimated by
\begin{eqnarray*}
&&C(T,g,L) \int_{-T}^T E(t,v) \d t + 2\int_{\Omega_T^-} \left| \partial_t v f
\right| \d t \d \nu \\
&&\leq C(T,g,L) \int_{-T}^T E(t,v) \d t + \| \partial_t v \|^2_{L^\infty
(]-T,T[\, ;~ L^2 (X))} + \| f \|^2_{L^1 (]-T,T[\, ;~ L^2 (X))} \\
&&\leq C(T,g,L) \int_{-T}^T E(t,v) \d t + \sup_{t\in ]-T,T[} E(t,v) + \| f
\|^2_{L^1 (]-T,T[\, ;~ L^2 (X))} \, ,
\end{eqnarray*}
where $C$ has the required continuity properties. Estimate (\ref{Ineq0H2}) then
gives (\ref{Ineq1H2}).
\item {\bf Proof of (\ref{Ineq2H2}).}
For the converse inequality, for $\min \varphi \leq t \leq T$, we
integrate (\ref{EnergyDensity}) on the domain $\Omega_t^+ = \left\{
  \varphi (x) \leq s \leq t \right\}$, i.e. the set of points of
$\tilde{X}$ situated in the future of $\Sigma$ and in the past of
$X_t$. Following H\"ormander, we put
\[ E_\varphi (t,v) = \int_{\varphi (x) \leq t} \left( \left(
    \partial_t v (t,x) \right)^2 + g^{\alpha \beta} (t,x)
  \partial_\alpha v (t,x) \partial_\beta v (t,x) +v(t,x)^2 \right) \d
\nu \, .\]
We obtain
\begin{eqnarray*}
0 &=& -\int_{\Sigma \cap \Omega_t^+} \left( \left( \partial_t v
  \right)^2 + g^{\alpha \beta} \partial_\alpha v \partial_\beta v +
  v^2 \right) \dsig + E_\varphi (t,v) \\
&&- \int_{\Sigma \cap \Omega_t^+} 2 g^{\alpha \beta} \partial_t v
\partial_\alpha \varphi \partial_\beta v \dsig \\
&&+ \int_{\Omega_t^+} \left( 2 \partial_t v L v + 2 \gamma^{-1}
  \partial_\alpha \left( \gamma g^{\alpha \beta} \right) \partial_t v
  \partial_\beta v \right. \\
&& \hspace{0.5in} \left. - \left( \partial_t g^{\alpha \beta} \right)
  \partial_\alpha v \partial_\beta v -2 v \partial_t v - 2 f\partial_t
  v \right) \d t\d
\nu \, .
\end{eqnarray*}
The first three terms give
\begin{gather*}
E_\varphi (t,v) - \left\| v_{|_\Sigma} \right\|^2_{H^1 (\Sigma \cap
  \Omega_t^+ )}
\end{gather*}
and the remainder can be estimated by
\[ C (T,g,L) \int_{-T}^t E_\varphi (s,v) \d s + \int_{\Omega^+_T} |f^2| \d s \d
\nu \, ,\]
where $C$ has the appropriate continuous dependence on $T$, $g$,
$g^{-1}$ and $L$. Gronwall's inequality entails for $\max \varphi
< t \leq T$
\[ E(t,v) \leq \tilde{C} (T,g,L) \left( \left\| v_{|_\Sigma} \right\|^2_{H^1
(\Sigma )} + \| f \|^2_{L^2 (\Omega_T )} \right) \, ,\]
where $\tilde{C}$ again has the required continuity properties. Eventually,
estimate (\ref{Ineq0H2}) gives (\ref{Ineq2H2}).
\end{enumerate}
This concludes the proof of proposition \ref{3EstimatesRegSol}. \qed

Let us now consider, for $k\in \N$, the solution $u_k$ of
(\ref{ApproxL1k}) associated with the initial data $u_0^k$, $u_1^k$
defined above. Using standard energy estimates of type (\ref{EnEst1}), we see
that for all $T>0$, $\| u_k \|_{{\cal F},T}$ is bounded uniformly in $k$. Now
consider the equation satisfied by $u_k -u_l$, for $k,l\in \N$~:
\begin{equation} \label{Equk-ul}
\partial_t^2 \left( u_k -u_l \right) - g^{\alpha \beta}
\partial_\alpha \partial_\beta \left( u_k -u_l \right)
+\tilde{L}_1^k \left( u_k -u_l \right) = \left( \tilde{L}_1^l -
  \tilde{L}_1^k \right) u_l \, .
\end{equation}
Let $T>0$, estimate (\ref{Ineq0H2}) and the hypotheses on $\tilde{L}^k_1$ give
the existence of a constant $C>0$, independent of $k,l$ and of $t\in
[-T,T]$ such that, for all $t \in [-T,T]$,
\[ E \left( t,u_k -u_l \right) \leq C \left( \left\| u_0^k -u_0^l
\right\|^2_{H^1 (X)} + \left\| u_1^k - u_1^l \right\|^2_{L^2 (X)} + \left\|
\left( \tilde{L}_1^l -
  \tilde{L}_1^k \right) u_l \right\|^2_{L^1 (]-T,T[\, ;~ L^2 (X))} \right) \, .
\]
Since
\[ \left( \tilde{L}_1^l -\tilde{L}_1^k \right) u_l = ( p_0^k - p_0^l )
\partial_t u_l + ( p^\alpha_k - p^\alpha_l ) \partial_\alpha u_l + (q_k-q_l)u_l
\, ,\]
using (\ref{ConvCoeffs1}), the boundedness of $\{ u_l \}_l$ in ${\cal C}^0
([-T,T] \, ;~ H^1 (X)) \cap {\cal C}^1 ([-T,T] \, ;~ L^2 (X))$, a Sobolev
embedding $H^1 (X) \hookrightarrow L^{p_1} (X)$ with $p_1>2$ and
(\ref{ConvCoeffs2}) for $p=2p_1/(p_1-2)$ we see that $\{ u_k \}_k$ converges
in ${\cal C} \left( [-T,T] \, ,~H^1 (X) \right) \cap {\cal C}^1 \left( [-T,T] \,
,~L^2 (X) \right)$. Allowing $k$ to tend to $+\infty$ in equation
(\ref{ApproxL1k}), we see that the limit of $u_k$ is the solution $u$ of
(\ref{WaveEq}) associated with the data $^t(u_0 ,u_1 )$ at $t=0$ (this uses the
convergence just
established as well as (\ref{ConvCoeffs1}), (\ref{ConvCoeffs2}) and
(\ref{ConvInData})). This convergence and inequality (\ref{Ineq1H2})
for $u_k -u_l$ then give that the restriction of $u_k$ to $\Sigma$
convergences strongly in $H^1 (\Sigma )$ (here again we need to use the
convergence of $\left( \tilde{L}_1^l -\tilde{L}_1^k \right) u_l$ towards $0$ in
$L^1 \left( ]-T,T[ \, ,~L^2 (X) \right)$). Moreover, using once again
the convergence of $u_k$ in ${\cal C}^0 ([-T,T] \, ;~ H^1 (X)) \cap
{\cal C}^1 ([-T,T] \, ;~ L^2 (X))$ and standard trace theorems, we see
that ${u_k}_{|_\Sigma}$ converges to $u_{|_\Sigma}$ in $L^2 (\Sigma)$. By
uniqueness it follows that ${u_k}_{|_\Sigma}$ converges to
$u_{|_\Sigma}$ in $H^1 (\Sigma)$. This entails inequalities
(\ref{EnEst2}) and (\ref{EnEst3}) for $u$.
\begin{remark} \label{RemProofCorol}
It is in this part of the proof that we require a bit more regularity on the
metric and the coefficients of the first order terms than
in H\"ormander's proposed setting, since we need the convergence of $\left(
\tilde{L}_1^l -\tilde{L}_1^k \right) u_l$ towards $0$ in $L^1 \left( ]-T,T[ \,
,~L^2 (X) \right)$ as $k,l\rightarrow +\infty$. Note that all we need to
guarantee this convergence is the convergence of $p^0_k$ towards $p^0$ and of
$p^\alpha_k$ towards $p^\alpha$ in $L^1_\mathrm{loc} (\R_t \, ;~{\cal C}^0 (X
))$ and of $q_k$ towards $q$ in $L^1_\mathrm{loc} (\R_t \, ;~L^p (X ))$ for all
$p<\infty$. This is true as soon as $g$, $p^0$ and $p^\alpha$ belong to
$L^\infty_\mathrm{loc} (\R_t \, ;~ {\cal C}^0 (X))$ and $q \in
L^\infty_\mathrm{loc} (\tilde{X})$. All the rest of the proof is valid for a
Lipschitz metric and coefficients of $L$ in $L^\infty_\mathrm{loc} (\tilde{X})$.
Remembering that $p^0$ and $p^\alpha$ contain first order derivatives of the
metric, this remark entails theorem \ref{TheorMinReg}.
\end{remark}

This shows that under hypothesis {\bf (H2)}, the operator $\Tsig$,
that to a solution $u$ of (\ref{WaveEq}) in $\cal E$
associates the trace of the solution $u$ on $\Sigma$ (well defined
since the solution is in ${\cal F} \hookrightarrow H^1_\mathrm{loc}
(\tilde{X})$), is a one-to-one bounded linear operator from $\cal E$
to $H^1 (\Sigma)$. It remains to establish the surjectivity.

$\bullet$ {\bf $\Tsig$ is surjective.} Let $v \in H^1 (\Sigma )$, we
prove that there exists $u \in {\cal E}$ such that $v= \Tsig u$. To do
so, we adopt the same regularization procedure as in the
second step of the proof of theorem \ref{Cauchy}, i.e. we consider $\{
\gk \}_k$ and $\{ \oplk \}_k$ defined by {\it (a)} and {\it (b)} and
satisfying (\ref{Convgk})-(\ref{UpLowBoundgk}). We introduce
for each $k$ a regularized equation to which we can apply theorem
\ref{thmHorm}. Equation (\ref{RegulEq}) will not do because we cannot
guarantee that $\Sigma$ is weakly spacelike for $\gk$. In order to
make up for this, all we need to do is slow down the propagation speed
for (\ref{RegulEq}). We consider a sequence $\{ \lambdak \}_k$,
$\lambdak \rightarrow 1$ as $k\rightarrow +\infty$, $0<\lambdak <1$,
such that
\[ \forall k \, , ~\lambdak \gk^{\alpha \beta} \left( x , \varphi
  (x) \right) \partial_\alpha \varphi (x) \partial_\beta \varphi (x)
  < 1 ~\mathrm{almost~everywhere~on~}X \, .\]
We have automatically that $\{ \lambdak^{-1} \gk \}_k$ satisfies
(\ref{Convgk}), (\ref{Boundgk}) and (\ref{UpLowBoundgk}). For each
$k$, we define the regularized equation
\begin{equation} \label{RegulEqGoursat}
\frac{\partial^2 v}{\partial t^2}- \lambdak \gamma^{-1}
  \frac{\partial}{\partial x^\alpha} \left( \gamma \left( \gk^{\alpha
  \beta}\right) \frac{\partial v}{\partial x^\beta} \right) + \oplk v
  =0 \, .
\end{equation}
That is to say, we have slowed down the propagation speed so that
$\Sigma$ is now totally spacelike for each equation
(\ref{RegulEqGoursat}) (i.e. for each $k$). We denote by $u_k$ the
unique solution of (\ref{RegulEqGoursat}) in $\cal F$ such that
$\left( u_k \right)_{|_\Sigma} = v$, $\left( \partial_t u_k
\right)_{|_\Sigma} = 0$ (the existence and uniqueness of such
solutions is given by theorem \ref{thmHorm}). For each $k$, using
theorem \ref{thmHorm}, we have an energy estimate (\ref{EnEst3}) for
solutions in $\cal F$ of equation (\ref{RegulEqGoursat}). Using the
properties of the regularized metric $\gk$ and operator $\oplk$, among
which the equivalence (uniform in $k$ and locally uniform in $t$)
between the energy (\ref{Energy}) induced by $g$ and that induced by
$\lambdak^{-1} \gk$, we obtain that $\{ u_k \}_k$ is bounded in ${\cal
  C} ([-T ,T] \, ,~H^1 (X) ) \cap {\cal C}^1 ([-T,T] \, , ~L^2 (X))$
for any $T > \max \{ |\min \varphi | \, ,~|\max \varphi |
\}$. The rest of the proof follows \cite{Ho} with elements of the
proof of theorem \ref{Cauchy} to deal with the regularized metric and
operator $\oplk$. Extracting a subsequence if necessary, we can
conclude that $u_k$ converges in the following spaces
\[ H^1 (]-T,T[ \times X) -w \, , ~H^s(]-T,T[ \times
X)~\mathrm{for~all}~s<1\, ,~L^\infty (]-T,T[\, ;~H^1 (X))-w-*\, , \]
towards a function $u$, $\partial_t u_k$ converges towards $\partial_t
u$ in $L^\infty (]-T,T[\, ;~L^2 (X))-w-*$. The convergences of $\gk$,
$\oplk$ and $u_k$ allow us to interpret, as in the existence part of
the proof of theorem \ref{Cauchy}, the convergence of each term of
equation (\ref{RegulEqGoursat}) in a common distribution space, hence
$u$ satisfies equation (\ref{WaveEq}) in the sense of
distributions. Following again the proof of theorem \ref{Cauchy}, we
show that $u$ belongs to $\tilde{\cal F}$ and therefore to $\cal
F$. Moreover, the strong convergence in $H^s (]-T,T[\times X)$ for all
$s<1$ entails the convergence in $L^2 (\Sigma)$ of the trace of $u_k$
on $\Sigma$ towards the trace of $u$ on $\Sigma$. Hence, $u_{|_\Sigma}
= v$. This concludes the proof of theorem \ref{Goursat}. \qed

\vspace{0.1in}
\begin{center}
{\Large {\bf Acknowledgements}}
\end{center}

The author would like to thank Alain Bachelot and Luc Robbiano for
helpful discussions while this work was in progress.

\end{document}